\documentclass[a4paper,10pt]{article}
\usepackage{amsmath} 

\usepackage{hyperref}
\usepackage{physics}
\usepackage{siunitx}
\usepackage{amssymb}
\usepackage{xcolor}
\usepackage{caption}
\usepackage{subcaption}
\usepackage{mathtools}
\usepackage{graphicx}
\usepackage{textcomp}
\usepackage{xcolor}

\newtheorem{remark}{\textbf{Remark}}

\newcommand{\tran}{^{\mkern-1.5mu\mathsf{T}}}

\newcommand{\mr}[1]{\mathrm{#1}}

\definecolor{dkgreen}{rgb}{0,0.4,0}

\usepackage[normalem]{ulem}

\title{\LARGE \bf
Controlled Synchronization of Coupled Pendulums by Koopman Model Predictive Control
}

\author{Loi Do$^1$, Milan Korda$^{1,2}$, and Zdeněk Hurák$^1$}

\begin{document}

\maketitle

\footnotetext[1]{Faculty of Electrical Engineering, Czech Technical University in Prague,
Technick\'a 2, CZ-16626 Prague, Czech Republic. {\tt \{doloi,hurak\}@fel.cvut.cz}}
\footnotetext[2]{CNRS; LAAS; 7 avenue du colonel Roche, F-31400 Toulouse; France. {\tt korda@laas.fr.}}

\begin{abstract}

We propose and experimentally demonstrate a feedback control method that allows synchronizing the motion of a chain of several coupled nonlinear oscillators actuated through one end of the chain. 
The chain considered in this work is a one-dimensional array of pendulums pivoting around a single axis and interacting with adjacent pendulums through torsion springs; the array is actuated using a single torque motor attached to one of the two boundary pendulums.
This represents a mechanical realization of the Frenkel-Kontorova model -- a spatially discrete version of a sine-Gordon equation describing (nonlinear) waves. 
The main challenges of controlling these systems are: high order (the number of pendulums can be high), nonlinear dynamics, and (as we set the problem here) only one actuator.
The presented problem of synchronization of motion is a special case of the problem of reference tracking, where all pendulums reach a common point or a trajectory.
In particular, we demonstrate synchronization to a stable equilibrium (all pendulums downward), unstable equilibrium (all pendulums upward), and a periodic orbit (all pendulums revolving).
We use the Koopman Model Predictive Control (KMPC) that constructs a linear predictor of the nonlinear system in a higher-dimensional lifted space and uses the predictor within a classical linear MPC, thereby maintaining low computational cost that allows for a real-time implementation, while taking into account the complex nonlinear dynamics.

\end{abstract}


\section{Introduction}

Studying the synchronization of coupled pendulums has a long history.
The foundations were already laid down in the 17th century by Christiaan Huygens, who described the synchronization of two pendulum clocks coupled through vibrations of the common base in~\cite{huygens_horologium_1673}.
Since then, the synchronization phenomena have been observed and studied in various systems of other coupled oscillators, from electrical circuits~\cite{van_der_pol_vii_1927},~\cite{huntoon_synchronization_1947} through biological systems~\cite{winfree_biological_1967},~\cite{strogatz_coupled_1993} to nanomechanical resonators~\cite{matheny_phase_2014}.

Adopting the terminology from~\cite{blekhman_self-synchronization_1997}, synchronization can be divided into \textit{self-synchronization} and \textit{controlled synchronization}.
The former describes a situation when synchronization between the subsystems occurs naturally without external inputs, induced directly by the character of the coupling.
This, for example, happens in many biological and physical systems, where the coupling is \textit{diffusive}, causing dissipation of the differences between the subsystems, thus eventually reaching a common state.
However, in some engineering applications, it is desired to directly control the synchronization, especially when the system does not self-synchronize or when one wishes to drive the system into a different synchronous state than the one reached through self-synchronization.
Thus, given an interconnected system, the task of controlled synchronization is to find a (feedback) control action that enables the states of an individual subsystem to approach a common point or a trajectory.

\begin{figure}[t]
        \centering      
        \includegraphics[width=0.7\textwidth]{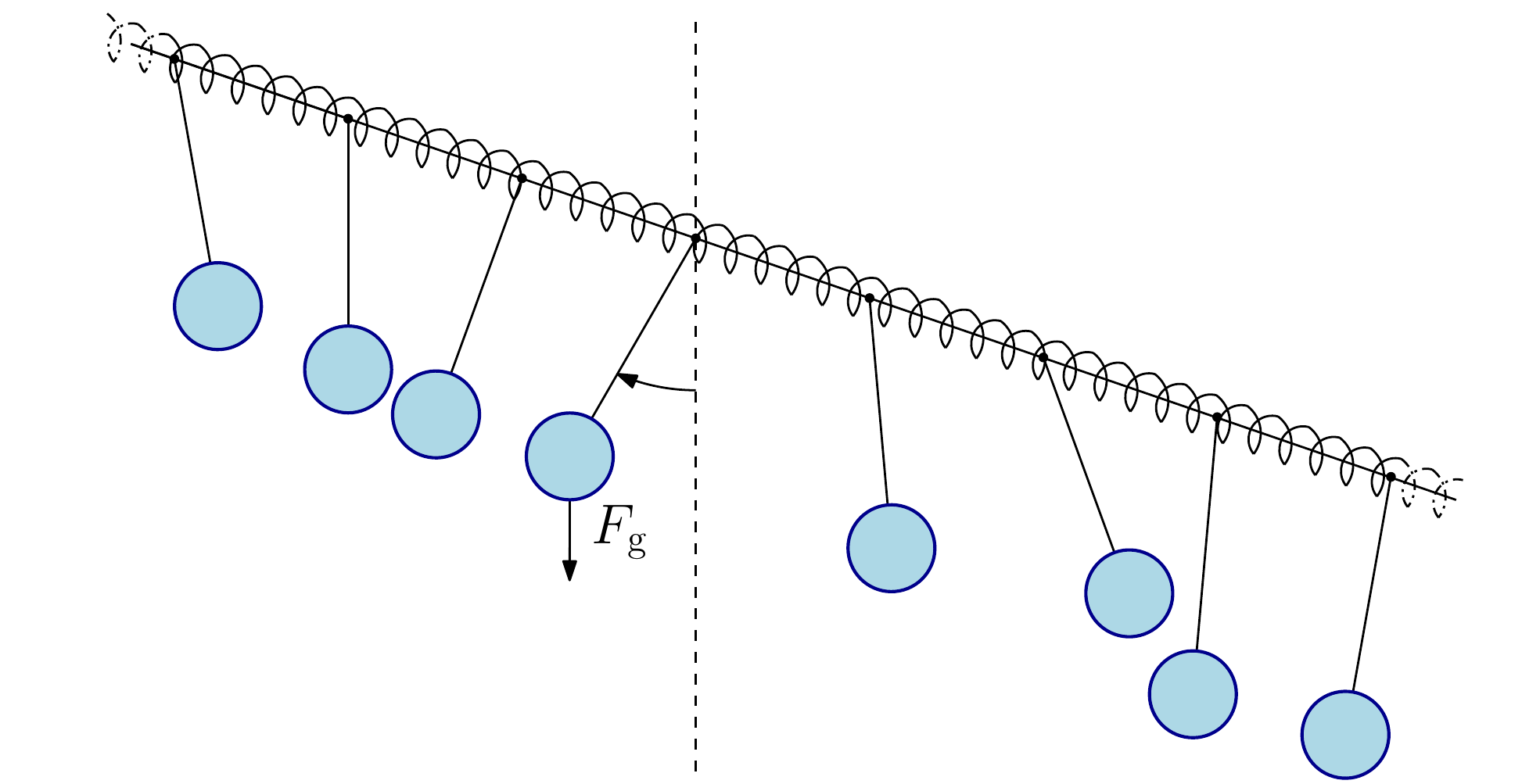}
        \caption{A series of coupled pendulums pivoting around a single axis}
        \label{fig:FK_chain_pendulums}
\end{figure}

This work deals with controlled synchronization in the system of coupled pendulums introduced in 1969 by A. C. Scott~\cite{scott_nonlinear_1969}, who constructed a series of pendulums pivoting around a single axis, where the adjacent pendulums were connected by torsion springs; see Fig.~\ref{fig:FK_chain_pendulums} for an illustration.
The author used the mechanical platform to demonstrate various wave solutions of the \textit{sine-Gordon} equation.
The model describing the system's dynamics was later coined as the \textit{Frenkel-Kontorova} (FK) model.

\paragraph{Motivation}

\begin{figure}[t]
        \centering      
        \includegraphics[width=0.7\textwidth]{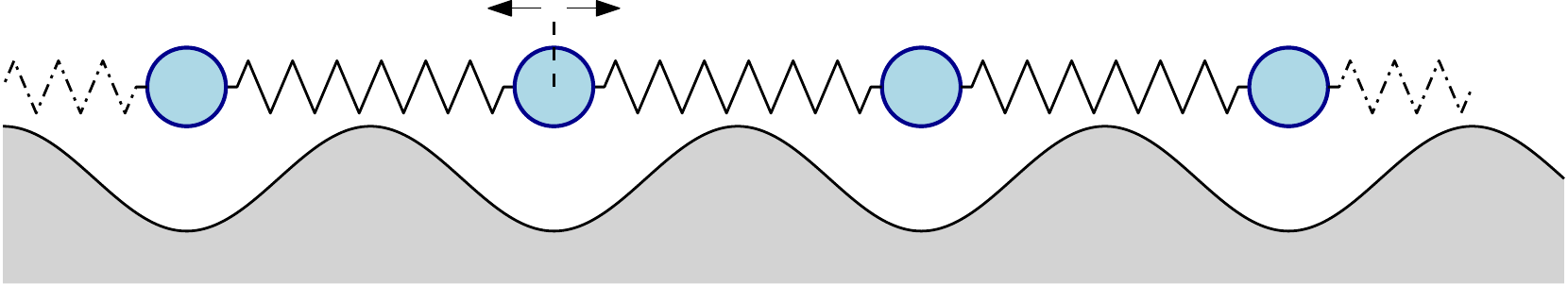}
        \caption{The FK model as an array of particles in a potential field}
        \label{fig:FK_chain_mass}
\end{figure}

The FK model was first presented in~\cite{frenkel_theory_1938} within a study of crystal dislocation, describing a motion of an infinite array of coupled identical particles in a spatially periodic potential field, see Fig.~\ref{fig:FK_chain_mass}.
Subsequently, the FK model has become one of the most fundamental models in physics, as it was found to describe a broad spectrum of various phenomena.
Apart from the crystal dislocation, it can be used to describe, for instance, the Josephson junction arrays, Bloch wall motion in magnetic domains, mechanical properties of an open-state of the DNA, or nanoscale friction.
Another significance of the FK model is that it represents a spatially discretized version of the \textit{sine-Gordon} (sG) equation, a completely integrable nonlinear partial differential equation.
Further details can be found in~\cite{braun_nonlinear_1998} or~\cite{cuevas-maraver_sine-gordon_2014}.

Over the years, many works have dealt with synchronizing coupled pendulums or motion control of the FK model.
Influenced by the work of Huygens, some researchers studied the controlled synchronization of two coupled pendulums with an additional objective of controlling the energy level in the system; see~\cite{kumon_controlled_2002},~\cite{pogromsky_controlled_2003}, and~\cite{fradkov_synchronization_2007}.

The works~\cite{thakur_driven_2007} and~\cite{cuevas-maraver_discrete_2009} have been dedicated to the analysis and experimental generation of particular solutions of the sine-Gordon equation, so-called \textit{Intrinsic Localized Modes}, also referred to as discrete \textit{breathers}.
The experimental platform used in the works was an array of pendulums -- the FK model -- where the whole frame was additionally driven in the horizontal direction by a motor. 
The authors showed that the array forms a stable \textit{discrete breather} when the platform's frame is driven with a suitable (open-loop) sinusoidal signal.

\paragraph{Main Contribution} 

In this paper, we formulate and solve three novel control problems on a finite array of $N$ pendulums modeled by the FK model.
In particular, we show, both in simulations and hardware experiments, the synchronization of the FK model into a stable equilibrium, an unstable equilibrium, and a periodic orbit.

Three factors make the considered control problems challenging. 
First, the system's dynamics is nonlinear and of high order. 
Second, we restrict the system's actuation only to \textit{boundary control}.
That is, the system of $N$ pendulums is controlled only by a single actuator, a torque motor attached to one of the two boundary pendulums in the array.
Therefore, the number of control inputs is significantly lower than the system's degrees of freedom.
Third, we consider a \textit{weak} coupling between the pendulums -- springs with relatively low stiffness, which makes the system more flexible and thus more difficult to control.

To deal with the challenges and solve the defined problems, we use the \textit{Koopman Model Predictive Control} (KMPC) introduced in~\cite{korda_linear_2018}. This method uses a higher-dimensional linear predictor of the  nonlinear dynamics, thereby  capturing the complex features of the nonlinear dynamics while allowing one to use a computationally cheap linear model predictive control (MPC) in closed-loop.

\subsection{Outline}

In Sec.~\ref{sec:mathematical_model}, we give a mathematical description of the FK model in the form of state-space equations and formulate the control problem to be solved in this paper.
We describe the method for solving the problem based on the Koopman Model Predictive Control, including a brief theoretical background in Sec.~\ref{sec:lin_predictor} and~\ref{sec:MPC}.
Then, we demonstrate the described control method through simulations and hardware experiments in Sec.~\ref{sec:control_sync}.
Finally, we discuss the results and outline the future work in Sec.~\ref{sec:future_work}.


\section{Model Description and Problem Statement}\label{sec:mathematical_model}

We consider the mechanical realization of the FK model in the form of an array of coupled and identical pendulums with a finite number of pendulums and with a single motor that applies torque directly to the first pendulum in the array.

\subsection{Mathematical Model}

Let $N$ be the number of pendulums in the array.
The equation describing the motion of the $i$-th pendulum, ${i = 1, \ldots, N}$, is
\begin{equation}\label{eq:main_model} 
\begin{split}
I\ddot{\varphi_i}       + mgl\sin(\varphi_i) + \gamma \dot{\varphi_i} - \frac{k}{2}\pdv{\varphi_i}\sum_{j=1}^{N-1}(\varphi_{j+1} - \varphi_j)^2 
                        - \frac{b}{2}\pdv{\dot{\varphi}_i} \sum_{j=1}^{N-1} (\dot{\varphi}_{j+1} - \dot{\varphi}_j)^2  =  M_i \;,
\end{split}
\end{equation}
where $g$ is the gravity constant, $k$ is the spring constant, and $I$, $m$, and $l$ are the moment of inertia, mass, and length of the pendulum, respectively.
The terms $M_i$ represent external torques. 
Since the system is driven only by one motor that applies the torque $M_1$ on the first pendulum in the array, we have $M_i = 0$ for $i = 2, 3, \ldots N$. 
Thus, the equations describing the motion of boundary pendulums can be written as
\begin{equation}
   \begin{aligned}
     I\ddot{\varphi_1}       + mgl\sin(\varphi_1) + \gamma \dot{\varphi_1} - k(\varphi_{2} - \varphi_1) 
                      - b(\dot{\varphi}_2 - \dot{\varphi}_1)  &=  M_1 \;, \\
     I\ddot{\varphi_N}       + mgl\sin(\varphi_N) + \gamma \dot{\varphi_N} - k(\varphi_{N-1} - \varphi_N) 
                      - b(\dot{\varphi}_{N-1} - \dot{\varphi}_i )  &= 0 \;,
   \end{aligned}
\end{equation}
and the equation for $i = 2, 3, \ldots N-1$ pendulums as
\begin{equation}
     I\ddot{\varphi_i}       + mgl\sin(\varphi_i) + \gamma \dot{\varphi_i} - k(\varphi_{i+1} - 2\varphi_i + \varphi_{i-1}) 
                       - b(\dot{\varphi}_{i+1} - 2\dot{\varphi}_i + \dot{\varphi}_{i-1})  =  0 \;.
\end{equation}
Numerical values of the physical parameters, that were identified from experiments with the hardware platform are in Tab.~\ref{tab:platform_parameters}.

\begin{table}[!tb] 
        \begin{center}
        \caption{Mechanical parameters of the FK model}\label{tab:platform_parameters}
        \begin{tabular}{lcc}
        Description                     & Symbol        & Value \\\hline
        Rod length                      & $l$           & $\SI{0.15}{\metre}$\\
        Spring constant                 & $k$      & $\SI{0.065}{\newton\per\metre}$ \\
        Pendulum's weight               & $m$           & $\SI{17}{\gram}$ \\
        Moment of inertia               & $I$           & $\SI{3.82e-4}{\kilo\gram\metre^2}$  \\
        \hline
        Relative dissipation coef.         & $b$           & $\SI{1.70e-3}{\newton\meter\second\per\radian}$  \\ 
        Absolute dissipation coef.         & $\gamma$      & $\SI{3.75e-4}{\newton\meter\second\per\radian}$
        \end{tabular}
        \end{center}
\end{table}


\subsection{FK model in the State-space Form}
For further development, it is convenient to rewrite the model~\eqref{eq:main_model} into a state-space form using the multi-agent system (MAS) formalism introduced in~\cite{fax_information_2004}.
The topological structure of the interactions between the pendulums can be described by an undirected path graph with a graph Laplacian ${L = (L_{ij}) \in \mathbb{R}^{N \times N}}$:
\begin{equation}\label{eq:Laplacian}
        L = 
        \begin{bmatrix}
        1     & -1    & 0     &     &\cdots& &  0     \\
        -1    & 2     & -1    &     &\cdots& &  0     \\
        0     &-1     & 2     &     &\cdots& &  0     \\
        \vdots&       &  & \ddots    &     &  &  \vdots\\
        0     &       &\cdots& &2    & -1    &  0     \\
        0     &       &\cdots& &-1   & 2     &  -1     \\
        0     &       &\cdots& &0    & -1    &  1   
        \end{bmatrix}\;.
\end{equation}
By defining the state of the $i$-th pendulum as ${x_i = [\varphi_i, \dot{\varphi}_i]\tran = [x_{i,1}, x_{i,2}]\tran}$, the state-space equations describing the dynamics can be written as
\begin{equation}\label{eq:FK_model_MAS}
  \dot{x}_i = f\left(x_i\right) - \frac{1}{I} G \left( K \sum_{j = 1}^N L_{ij} x_j - M_i \right)\;,
\end{equation}
where $f(x_i)$ is the uncoupled dynamics of a single pendulum
\begin{equation}\label{eq:drift_dynamics_MAS}
  f(x_i) =
  \begin{bmatrix}
    x_{i,2} \\
    -\dfrac{mgl}{I}\sin(x_{i,1}) -\dfrac{\gamma}{I} x_{i,2} 
  \end{bmatrix}\;,
\end{equation}
and ${G = [0, 1]\tran}$, ${K = [k, b]}$, which represent the torsion coupling through springs with dissipation.
Additionally, let 
\begin{equation}\label{eq:MAS_stacked_vectors}
        \begin{split}
                u &= M_1 \;, \\
                x &= [x_1\tran, \ldots, x_N\tran]\tran           \;,    \\
                F(x) &= [f(x_1)\tran, \ldots, f(x_N)\tran]\tran  \;,    \\
                d &= [1, 0, \ldots, 0]\tran \in  \mathbb{R}^{N} \;, \\
                D &= \text{diag}(d) \;,
        \end{split}
\end{equation}
where $\mr{diag}(d)$ is a diagonal matrix with element of the vector $d$ on the main diagonal.
Using~\eqref{eq:MAS_stacked_vectors}, we can describe the system in a total state-space form
\begin{equation}\label{eq:matrix_form_model_main}
        \begin{split}
                \dot{x} &= F(x) - \left( \frac{1}{I} (L+D) \otimes G K \right) x  + \frac{1}{I}(d \otimes G) u \;, \\
                y &= x \;, 
        \end{split}
\end{equation}
where $\otimes$ denotes the Kronecker product and $y$ is the (measured) output of the system.
Thus, we assume that all states are measured.


\subsection{Problem Statement}\label{sec:problem_statement}

We consider the problem of synchronizing the individual states $x_i$ in the system~\eqref{eq:matrix_form_model_main}.
In particular, we denote $x^\star(t)$ to be the synchronization trajectory and let $x^\star(t)$ be a solution to a single, \textit{isolated} (uncoupled) pendulum, satisfying
\begin{equation}\label{eq:virtual_leader}
    \dot{x}^\star(t) = f(x^\star) \;,
\end{equation}
where $f(x^\star)$ is the drift dynamics~\eqref{eq:drift_dynamics_MAS}.
The solution $x^\star(t)$ may be an equilibrium point or a periodic orbit.
In the framework of multi-agent systems, the system~\eqref{eq:virtual_leader} is referred to as a \textit{virtual leader}. 

The goal is to find a closed-loop controller $u = M_1(x)$ such that all pendulums in the chain synchronizes with the trajectory $x^\star(t)$, in the sense that
\begin{equation}\label{eq:main_problem_statement}
        \lim_{t \rightarrow \infty} | x_i(t) - x^\star(t) |  < \varepsilon \;, \quad i = 1, 2, \ldots, N \;,
\end{equation}
from any initial conditions, and $\varepsilon > 0$ is a small constant value.
By introducing the constant $\varepsilon$, we allow non-zero synchronization error, thus requiring only \textit{approximated synchronization}.

To achieve synchronization on a nontrivial (non-zero) trajectory in a system of interconnected, nonlinear system, the reference trajectory $x^\star(t)$ has to be necessarily the solution to the drift dynamics, see~\cite{wieland_static_2010}.
The basic explanation is the following.
In the FK model, a spring connecting two adjacent pendulums create a \textit{diffusive} coupling, i.e.,  the coupling naturally dissipates the state differences and drives the individual states $x_i$ to a common state.
If the differences dissipate, the system synchronizes and both the control $u$ and coupling vanish, so the individual subsystems behave as \textit{uncoupled}, thus, driven only by their drift dynamics.
For more rigorous description and further details, see, for instance,~\cite{delellis_quad_2011}, ~\cite{yu_synchronization_2013}, or~\cite{hengster-movric_structured_2016}, and references in them.

\begin{remark}
    The defined problem is a special case of the trajectory tracking.
    However, instead of considering any reference trajectory for the total system's state $x \in \mathbb{R}^{2N}$, we restrict ourselves to a particular trajectories where $x(t) = [x^\star(t), x^\star(t), \ldots, x^\star(t)]\tran$ for $t \rightarrow \infty$. 
\end{remark}

To solve the defined problem~\eqref{eq:main_problem_statement}, we use the Koopman Model Predictive Control introduced in~\cite{korda_linear_2018}.
That is, we use the linear (discrete time) MPC with a tailored linear predictor of the system's nonlinear dynamics based on the Koopman operator theory.
We construct the linear predictor using the data-driven method, the Extended Dynamic Mode Decomposition (EDMD) algorithm.


\section{Linear Predictor based on the Koopman Operator Theory}\label{sec:lin_predictor}

The main idea behind obtaining accurate predictions of nonlinear dynamics using a linear model is to \textit{lift} the nonlinear dynamics into a higher dimensional space where the evolution of the lifted state is approximately linear.

We seek the predictor for the system in the form of a discrete-time linear system
\begin{equation}\label{eq:linear_predictor}
        \begin{split}
                z_{k+1} &= Az_k + Bu_k \;,      \quad z \in \mathbb{R}^{n}\;, u \in \mathbb{R} \;,\\
                \hat{y}_k &= Cz_k \;,           \quad  \hat{y} \in \mathbb{R}^{2N} \;,
        \end{split}
\end{equation}
starting from an initial condition $z_0=\boldsymbol{\psi}(x_0)$, where $\boldsymbol{\psi}:\mathbb{R}^{2N}\to \mathbb{R}^{n}$, with $n\ge 2N$, is a user-specified nonlinear (lifting) mapping, and $\hat{y}_k$ is the prediction of the nonlinear system's output $y_k$ at time $k$.
The main challenge is to find the matrices $A$, $B$, $C$, and $\boldsymbol\psi$ that represent the system's behavior well.

\subsection{Koopman Operator}

To rigorously justify the construction of the linear predictor, we shortly describe the \textit{Koopman operator} approach for analysis of dynamical systems. This approach dates back to the seminal works of Koopman and von Neumann~\cite{Koopman1931, koopman1932dynamical} in the 1930s, with resurgence of interest from the mid 2000s starting with the works~\cite{mezic_spectral_2005,Mezic_control_suggestion}.
Consider an uncontrolled discrete-time dynamical system
\begin{equation}\label{eq:disc_non-linear_sys_Koopman}
        x_{k+1} = \mathcal{T}(x_k) \;, \quad x_k \in \mathcal{M} \;,
\end{equation}
where $\mathcal{T}(x_k)$ is, in general, a non-linear transition mapping, and $\mathcal{M}$ is the state space.
Note that the presented ideas can be formulated also for systems with external inputs acting as control; this was first formalized and used in conjunction with model predictive control in \cite{korda_linear_2018}.

Instead of analyzing the mapping $\mathcal{T}$, we investigate how functions of the states, so-called \textit{observables} evolve along the \textit{flow} of the system.
Formally, an observable is a function ${\psi:\mathcal{M} \rightarrow \mathbb{R}}$ belonging to a suitable, typically infinite-dimensional, space $\mathcal{F}$.
The \textit{Koopman operator} ${\mathcal{K}:\mathcal{F}\rightarrow\mathcal{F}}$ is then defined as 
\begin{equation}
        ( \mathcal{K} \psi )(x_k) = \psi (T(x_k)) = \psi(x_{k+1}) \;,
\end{equation}
Thus, the (discrete-time) Koopman operator $\mathcal{K}$ advances the observables $\psi$ from the time $k$ to the next time step $k+1$.
Note that the space $\mathcal{F}$ is invariant under the action of $\mathcal{K}$.

There are two key properties of the Koopman operator.
First, it fully captures the behavior of the original non-linear dynamical system, provided that the space of observables $\mathcal{F}$ contains the coordinate identity mappings $x\mapsto x_i$.
Second, the operator is linear since for any two observables $\psi_1$ and $\psi_2$, and scalar values $a_1$ and $a_2$, it holds
\begin{equation}
        \mathcal{K} (a_1\psi_1 + a_2\psi_2) = a_1 \mathcal{K}(\psi_1) + a_2 \mathcal{K}(\psi_2) \;.
\end{equation}
Therefore, using the Koopman operator, we have converted the analysis of the finite-dimensional non-linear system into the analysis of the infinite-dimensional linear operator. 
For some systems, one can find an invariant subspace of $\mathcal{F}$ that is spanned by a finite number of observables while the observables exactly capture the systems' dynamics (i.e., they contain the coordinate identity mappings in their span).
When the basis of this subspace is fixed, the operator $\mathcal{K}$ restricted to this subspace can be represented by a matrix of a finite size. For further details about the Koopman operator theory, we refer the reader to the two surveys~\cite{budisic_applied_2012,brunton2021modern}.

\subsection{Extended DMD}

A problem with using the Koopman operator is that for most non-linear systems, including the FK model~\eqref{eq:matrix_form_model_main}, finding an invariant subspace of $\mathcal{K}$ containing the observables we wish to predict is intractable.
However, we can look for its finite approximation, which then induces a linear predictor in the form~\eqref{eq:linear_predictor}.
With a suitably selected set of observables, the linear predictor can capture the original system's dynamics with high accuracy; see~\cite{korda_linear_2018} for the underlying theory of how to achieve this with control.
To find the predictor, we use the EDMD algorithm introduced in~\cite{williams_datadriven_2015} with an extension for controlled systems proposed in~\cite{korda_linear_2018}. The convergence of EDMD was proven in~\cite{korda2018convergence} and analyzed quantitatively in \cite{zhang2023quantitative,nuske2023finite}. A more refined way to construct the predictors exploiting the system dynamics and the state-space geometry was proposed in~\cite{korda2020optimal}.

We first select $n$ observables (the lifting functions) $\psi_i$, and form the predictor's state from the output $y$ of the system~\eqref{eq:matrix_form_model_main} as
\begin{equation}\label{eq:n_observables_general}
        z(y) = [\psi_1(y), \psi_2(y), \ldots, \psi_n(y)]\tran\;.
\end{equation}
Next, we gather measurements $y_i$ and $y_i^+$ from the system and form a set of data 
\begin{equation}\label{eq:EDMD_matrices}
\begin{split}
    X           & = [y_1, y_2, \ldots, y_{N_\mr{d}}]\;, \\
    X_\mr{lift}  & = [z(y_1), z(y_2), \ldots, z(y_{N_\mr{d}})]\;, \\
    Y_\mr{lift}  &= [z(y_1^+), z(y_2^+), \ldots, z(y_{N_\mr{d}}^+)]\;, \\
    U           &= [u_1, u_2, \ldots, u_{N_\mr{d}}]\;. 
\end{split}
\end{equation}
where the measurements satisfy a relation $y_i^+ = \mathcal{T}(y_i,u_i)$ that we wish to capture with the linear predictor, and $\mathcal{T}(y_i,u_i)$ is an input-dependent transition mapping.
Note, that the outputs in the set of data~\eqref{eq:EDMD_matrices} need not be temporally ordered, i.e., the data could be gathered from multiple trajectories of the system~\eqref{eq:matrix_form_model_main}.

The linear predictor~\eqref{eq:linear_predictor} can be then identified by solving
\begin{subequations}\label{eq:DMD_problem_def}
        \begin{equation}
                \min_{A,B}      \| Y_\mr{lift} - AX_\mr{lift} - BU \|_\mr{F} \;,
        \end{equation}
        \begin{equation}\label{eq:DMD_problem_def_C_mat}
                \min_{C}        \| X - C X_\mr{lift} \|_\mr{F}  \;,
        \end{equation}
\end{subequations}
where $\| . \|_\mr{F}$ denotes a Frobenius norm of a matrix.
The solution to~\eqref{eq:DMD_problem_def} can be analytically obtained from
\begin{equation}
        \begin{bmatrix}
                A & B \\
                C & 0 
        \end{bmatrix}
        =
        \begin{bmatrix}
                Y_\mr{lift} \\
                X
        \end{bmatrix}        
        \begin{bmatrix}
                X_\mr{lift} \\
                U
        \end{bmatrix}\tran 
        \left(
        \begin{bmatrix}
                X_\mr{lift} \\
                U
        \end{bmatrix}
        \begin{bmatrix}
                X_\mr{lift} \\
                U
        \end{bmatrix}
        \tran 
        \right)^\dagger     \;,
\end{equation}
where $(.)^\dagger$ denotes the Moore-Penrose pseudoinverse of a matrix.
Note, that when all original states $x$ of the non-linear system are incorporated in the set of observables, the matrix $C$ can be directly constructed by selecting the corresponding observables from~\eqref{eq:n_observables_general}.


\section{Koopman Model Predictive Control}\label{sec:MPC}

The discrete-time linear Model Predictive Control (MPC) is a closed-loop technique that finds the control inputs by solving an optimization problem that minimizes a user-defined cost function $J$ reflecting the control goal.
For a standard linear MPC, the optimization problem is a convex quadratic program where the model of the system is used to predict the system's behavior on a \textit{prediction horizon} $N_\mathrm{p}$, allowing one to find the optimal control sequence $u_k$ on this horizon.
Only the first input of the optimal sequence is then applied to the system, and the optimization problem is resolved in the next time step; this is often referred to as a receding horizon approach.
The KMPC \cite{korda_linear_2018} extends this technique by using the linear predictor~\ref{eq:linear_predictor} constructed by the Koopman operator theory to predict the system's nonlinear behavior instead of using a local linearization of the system's model.

Specifically, to formulate the MPC for reference tracking, let $e_k = r_k - \hat{y}_k$ be a tracking error and $r_k$ a given reference signal having the same dimension as ${\hat{y}_k = Cz_k}$.
The optimization problem solved at every time step $t$ is
\begin{equation}\label{eq:sparse_MPC_QP}
        \begin{aligned}
                \min_{u_k, e_k}         \quad   & J\left( \{e_k\}_{k = 0}^{N_\mr{p}}, \{u_k\}_{k = 0}^{N_\mr{p}-1} \right) \;, \\
                \textrm{subject to}     \quad   & z_{k+1} = Az_k + Bu_k         \;, \quad       k = 0, 1, \ldots, N_\mathrm{p} - 1 \;,\\
                                                & e_k = r_k - Cz_k                 \;,\\
                                                & z_\mr{min} \leq z_k \leq z_\mr{max}   \;,\\
                                                & u_\mr{min} \leq u_k \leq u_\mr{max}   \;,\\
                \text{parameters}       \quad   & z_0 = \psi(y_t)                       \;,\\
                                                & r_k = \text{given}  \;, \quad k = 0, 1, \ldots, N_\mathrm{p} \;,\\
        \end{aligned}
\end{equation}
where $y_t$ is the measured output from the original system at time $t$ and the terms with subscripts '$\mr{min}$' and '$\mr{max}$' define the bounds on the optimization variables.
The cost function $J(.)$ has a standard quadratic form for MPC tracking 
\begin{equation}\label{eq:MPC_cost_function}
        \begin{split}
                J(e_k,u_k) &= e_{N_\mr{p}}\tran Q_N e_{N_\mr{p}} + \sum_{k = 0}^{N_\mr{p} - 1} 
                \left[
                        e_{k}\tran Q e_{k} + u_{k}\tran R u_{k}
                \right],
        \end{split}
\end{equation}
where $Q \succeq 0$, $Q_N \succeq 0$, and $R \succeq 0$ are cost matrices. 
Importantly, the optimization problem~\eqref{eq:MPC_cost_function} is a convex quadratic program that can be solved by high-performance software tailored to the MPC structure encountered here. 
We used the OSQP solver, see~\cite{stellato_osqp_2020}.

\begin{remark}
        The optimization problem could also be reformulated to penalize the increments of inputs $\Delta u_k$, instead of the inputs $u_k$. 
        This would allow to not penalize non-zero, constant input in the steady-state.
        However, this is not required in our case.
\end{remark}

\subsection{Dense Formulation}
To reduce the computational complexity of the formulated optimization problem~\eqref{eq:sparse_MPC_QP}, we use the \textit{dense} formulation of~\eqref{eq:sparse_MPC_QP}.
In particular, using the state-space equations~\eqref{eq:linear_predictor}, we can explicitly express the outputs~$\hat{y}_k$ as functions of~$u_k$.
Thus, the complexity of the optimization problem in the dense formulation is given only by the number of inputs instead of depending on the size of the lifted state~$z$.

In particular, we rewrite the optimization problem into a form
\begin{equation}\label{eq:MPC_formulation_dense}
         \begin{aligned}
                 \min_{U}          \quad &     \frac{1}{2} U\tran H U +  z_0\tran G U \;,\\
                 \textrm{s.t.}     \quad &     b_\mr{min} \leq \bar{{A}} z_0 + \bar{{B}} U \leq b_\mr{max} \;, \\
                 \mathrm{parameter}\quad &     z_0 = \psi(y_t) \;,
         \end{aligned}
\end{equation}
where the $U = [u_0\tran, \ldots, u_{N_\mathrm{p} - 1 }]\tran$. The matrices $H, G, {\bar{A}}, {\bar{B}}$, and vectors $b_\mr{min}$ and $b_\mr{min}$ are listed in Appendix~\ref{sec:Appdx:dense_MPC}.


\section{Controlled Synchronization in the FK Model}\label{sec:control_sync}

We now explore how the presented methodology can be used to synchronize the pendulums in the FK model.
In particular, we show synchronization into a stable equilibrium $x_\mathrm{stab}$, an unstable equilibrium $x_\mathrm{unst}$, and a periodic orbit~$x_\mathrm{per}$.
We assume, that the springs between the pendulums provide only a weak coupling, so the nonlinear drift dynamics of a single pendulum is dominant.  
This makes the problem more challenging, compared to the case with strong coupling, as the synchronization could not be reached by directly controlling only one pendulum, but the states of all pendulums and their complex dynamics need to be included exploited in the control.

\begin{figure}[t]
        \centering      
        \includegraphics[width=0.9\textwidth]{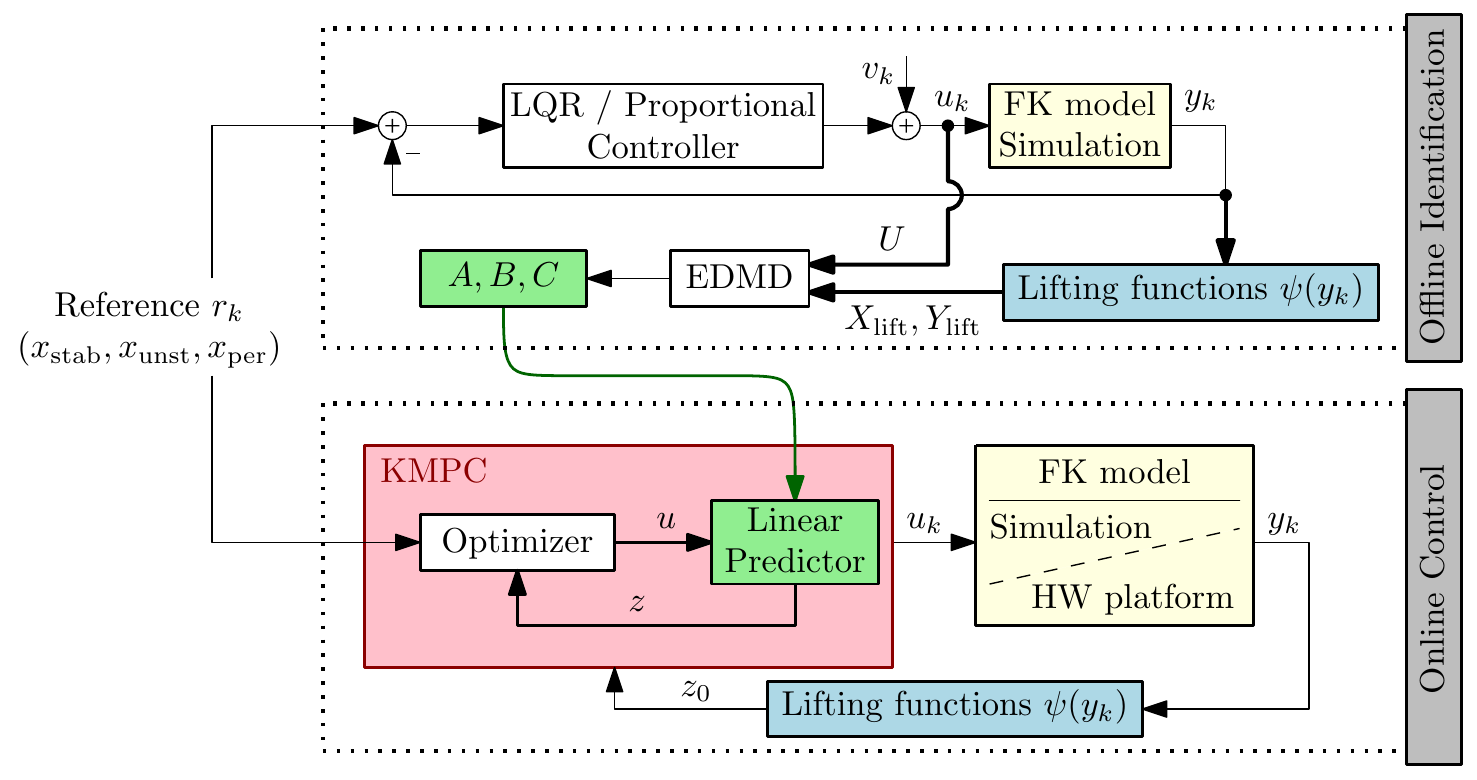}
        \caption{Schematic representation of the KMPC framework with closed-loop identification}
        \label{fig:big_picture_diagram}
\end{figure}

The schematic representation of the KMPC framework is in Fig.~\ref{fig:big_picture_diagram}.
We start by constructing the linear predictor for the FK model.
This step involves gathering data from the system and selecting a suitable set of observables.
We identify a different predictor for each task, directly tailored to the particular reference signal $r_k$.
Lastly, we run the KMPC algorithm as described in Sec.~\ref{sec:MPC}.

For the KMPC in all tasks, we do not put any constraints on the state of the predictor~\eqref{eq:linear_predictor}, but we constrained the input, both in simulation and on the hardware platform, to be in the interval ${-0.1 \leq u_k < 0.1}$.
The motivation for constraining the input is to prevent extensive torques from being applied to the springs, thus preventing their non-elastic deformation and potential damage.
Both the simulations and experiments were done with a sampling time $\Delta t = \SI{5}{\milli\second}$, given by hardware limitations.
The prediction horizon for the simulations was set to $N_\mathrm{p} = 50$, while for the experiments, we set the prediction horizon to $N_\mathrm{p} = 20$ due to hardware limitations.

\subsection{Identification Data}\label{sec:EDMD_identification}

A challenging step in the EDMD algorithm is to get appropriate trajectories for the system identification, i.e., to fill the matrices~\eqref{eq:EDMD_matrices} so the model accurately predicts the system's dynamics along the desired trajectories.

One option is to use common input signals for identification, such as square wave or chirp signals or their linear combination based on the a priori assumption of the closed-loop input.
However, the assumptions of closed-loop $u_k$ might not be correct and could result in model fitted to incorrect trajectories.  
The second option, applied in this work, is to use \textit{closed-loop} identification, where the identified system is controlled by a feedback controller rather than driven by an open-loop signal (see, e.g., \cite{forssell1999closed,hjalmarsson2005experiment}).

For the simulations, we first design an admissible feedback controller using techniques that do not require accurate predictions of the system.
This controller is then used to acquire the data from closed-loop simulations.
We describe the particular design of the controllers in Sec.~\ref{sec:Simulations}.
For the hardware experiments, we acquire the data from closed-loop simulations of the mathematical model instead of gathering data directly on the hardware platform.
Although there is an inevitable mismatch between the mathematical model and the hardware platform, this approach allows us to conveniently generate many different trajectories, i.e., with different initial conditions or random perturbations of the input.

\subsection{Observables for the FK Model}

The selection of suitable observables is the central challenge in the EDMD methodology, especially when the system's dynamics is not known.
In our case, the situation is simpler.
Not only is the system's dynamics known, but we can also exploit the structure of the FK model.

In particular, observing the FK model in the form~\eqref{eq:matrix_form_model_main}, the system consists of linearly coupled, identical non-linear subsystems -- the pendulums -- with the drift dynamics~\eqref{eq:drift_dynamics_MAS}.
For every pendulum in the FK model, we choose $n=6$ observables  
\begin{equation}\label{eq:set_observables_1}
        z_i = [\varphi_i, \dot{\varphi}_i, \sin(\varphi_i), \cos(\varphi_i), \dot{\varphi}_i\sin(\varphi_i), \dot{\varphi}_i\cos(\varphi_i)]\tran \;,
\end{equation}
so the total set of observables for the FK model with $N$ pendulums is
\begin{equation}\label{eq:set_observables_all}
        z = [z_1\tran, z_2\tran, \ldots, z_N\tran]\tran \in \mathbb{R}^{6N} \;.
\end{equation}
The rationale behind selecting~\eqref{eq:set_observables_1} as the observables is the following.
The linear part of the drift dynamics and linear coupling can be directly captured by including the pendulums' states $\varphi_i, \dot{\varphi}_i$ into the predictor's state.
By including the nonlinear observables, we wish to capture the time-evolution of the nonlinear term $\sin(\varphi_i)$ along the system's trajectory.
This selection of observables is also supported by the fact that as the FK-model consists of $N$-coupled oscillators, the system's trajectories are confined on an $N$-torus, which is parameterized by sine and cosine terms.

\subsection{Simulations}\label{sec:Simulations}


\subsubsection{Stable Equilibrium}\label{sec:task_equilibrium_stab}
One can easily check that the stable equilibrium of the FK model~\eqref{eq:matrix_form_model_main} is the origin
\begin{equation}
        x_\mathrm{stab} =  [0, 0, \ldots, 0]\tran \;,
\end{equation}
The goal is to satisfy~\eqref{eq:main_problem_statement} with ${x^\star = x_\mathrm{stab}}$ for any initial conditions. 
Note, that when the initial conditions are close to the $x_\mathrm{stab}$, the system naturally synchronizes even without external control.
Nevertheless, the applied control can  significantly speed up the convergence to the synchrony.
The practical motivation for this task is a control of vibrations -- attenuation of oscillations in flexible structures.


To solve the goal within the MPC framework, we selected the reference signal for the total system as
\begin{equation}
        r_k = [0, 0, \ldots, 0] \in \mathbb{R}^{2N}\;,
\end{equation}
driving both the angles and speeds of all pendulums to zero.

For the predictor's identification, we gather the closed-loop trajectories by designing Linear Quadratic Regulator (LQR) that minimizes the cost functional
\begin{equation}\label{eq:LQR_functional}
        \mathcal{J} = \int_0^{\infty} \left(  x\tran(t) Q x(t) + u(t)\tran R u(t) \right) dt \;,
\end{equation}
subjected to the dynamics of a linear approximation of the system~\eqref{eq:matrix_form_model_main}.
Specifically, the linear approximation around the $x_\mathrm{stab}$ yields a linear system
\begin{equation}
        \Delta \dot{x} = \tilde{A} \Delta x + \tilde{B} \Delta u \;, 
\end{equation}
with the system's matrices
\begin{equation}
        \tilde{A} =          
        \left(I_N
                \otimes
                \begin{bmatrix}
                        0 & 1 \\
                        -\frac{mgl}{I} & -\frac{\gamma}{I} 
                \end{bmatrix}
        \right) - \left( \frac{1}{I} (L+D) \otimes G K \right) \;, 
        \tilde{B} = \frac{1}{I}(d \otimes G) \;,
\end{equation}
where $I_N$ denotes an identity matrix of size $N$.
The LQR feedback input is given by 
\begin{equation}\label{eq:sync_ex1_lqr_ctrl}
        u_{k}^\mathrm{LQR} = -R^{-1}\tilde{B}\tran S x_k \;,        
\end{equation}
where $S$ is a solution to the continuous time algebraic Riccati equation.
The cost matrices in~\eqref{eq:LQR_functional} were set to
\begin{equation}
    \begin{split}
        Q &= \text{blkdiag} 
            \left( 
            \left\{
                \begin{bmatrix}
                    1000 & 0 \\
                    0 & 0.01
                \end{bmatrix}
            \right\}_{i=1}^N
            \right) \;, \\  
        R &= 0.1 \;,
    \end{split}
\end{equation}
where $\text{blkdiag}(.)$ is a block diagonal matrix of size $2N$.
To get richer trajectories, we additionally add random perturbations $v_k$ to~\eqref{eq:sync_ex1_lqr_ctrl}, so the resulted control action for identification is
\begin{equation}\label{eq:sync_ex1_lqr_ctrl_with_noise}
        u_k = u_k^\mathrm{LQR} + v_k \;, \quad v_k \sim \mathcal{N}(0, 0.1) \;,
\end{equation}
where $\mathcal{N}(\mu, \sigma^2)$ denotes normal distribution with mean $\mu$ and variance $\sigma^2$.
We then used the resulted control~\eqref{eq:sync_ex1_lqr_ctrl_with_noise} to simulate the system~\eqref{eq:matrix_form_model_main} starting from random initial conditions to gather $N_\mathrm{traj} = 200$ trajectories over 200 sampling periods to identify the predictor via EDMD algorithm described in the Sec.~\ref{sec:EDMD_identification}.
Each trajectory was $\SI{5}{\second}$ long and sampled with $\Delta t$.

To get the control using the MPC, we set the cost matrices in~\eqref{eq:MPC_cost_function} as
\begin{equation}\label{eq:Q_R_matrices_equilibrium_stab}
        \begin{split}
                Q &= \text{blkdiag} \left( \{Q_i\}^N_{i = 1} \right) \;, \\
                Q_N &= Q \;, \\
                R & = 0.1
                \;.
        \end{split}
\end{equation}
We set the blocks $Q_i$ as
\begin{equation}\label{eq:Q_matrix_increasing}
        Q_i = 
        \begin{bmatrix}
                10 i^3        & 0 \\
                0               & 0.01
        \end{bmatrix} \;.
\end{equation}
This choice allow us to incorporate the structure of the system into the optimization problem.
In particular, with increasing distance from the actuator (attached to the first pendulum ${i=1}$), the controllability of pendulums' states decreases.
Thus, the intention is to stabilize the array from the last to the first pendulum.

\begin{figure}[t]
        \centering      
        \includegraphics[width=0.7\textwidth]{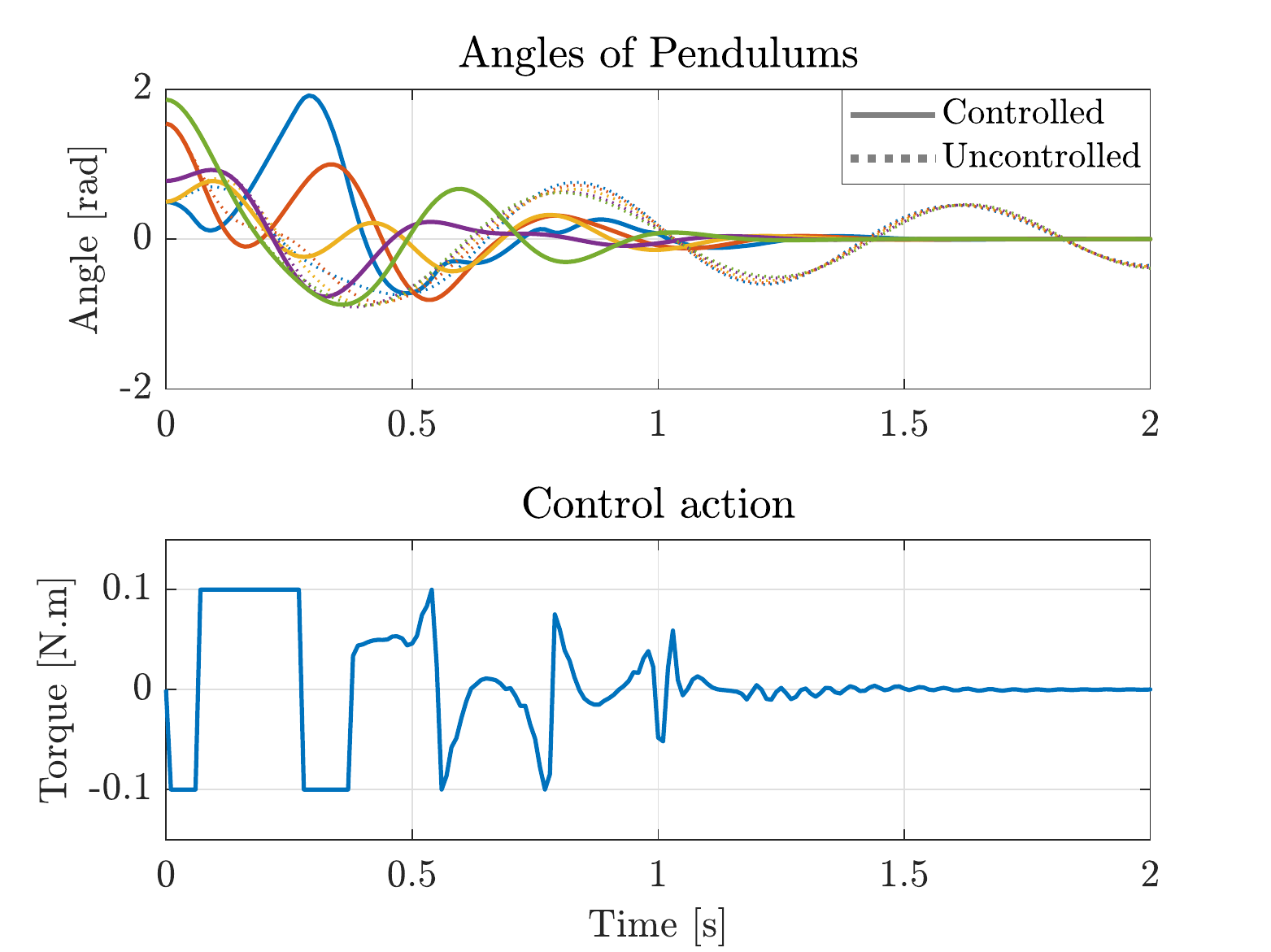}
        \caption{Controlled synchronization to the stable equilibrium $x_\mathrm{stab}$ compared with uncontrolled dynamics}
        \label{fig:sim_ctrl_5_pend_stable_equi}
\end{figure}

The resulting control is depicted in Fig.~\ref{fig:sim_ctrl_5_pend_stable_equi} with plotted angles of pendulums and the control action.
When compared to the uncontrolled response to the initial conditions (dashed lines), the controlled system converges to the $x_\mathrm{stab}$ significantly faster.

\subsubsection{Unstable Equilibrium -- Swing-up}\label{sec:task_swing_up}

The FK model also has an infinite number of unstable equilibria. 
We again focus only on synchronization to a single equilibrium.
It is straightforward to check, that the unstable equilibrium of the FK model~\eqref{eq:main_model} is
\begin{equation}\label{eq:unstable_equilibrium}
        x_\mathrm{unst} = [0, \pi, 0, \pi, \ldots, 0, \pi]\tran \in \mathbb{R}^{2N} \;.
\end{equation}
Reaching the unstable equilibrium~\eqref{eq:unstable_equilibrium} from the origin corresponds to the classical task in control theory -- the \textit{Swing-up}.
Many researchers have studied the swing-up task in pendulum-like systems as it nicely illustrates several concepts from the field of non-linear control (see, e.g.,~\cite{spong_swing_1995},~\cite{astrom_swinging_2000}).
The task is to satisfy~\eqref{eq:main_problem_statement} with $x^\star = x_\mathrm{unst}$ for the system starting from the origin.
Similar to the swing-up of multiple-link pendulums, the main challenge in our setup is that the system is underactuated, i.e. the number of degrees-of-freedom in the system is higher than number of inputs.

We select the reference trajectory $r_k$ for the total system as
\begin{equation}
        r_k = [0, \pi, 0, \pi, \ldots, 0, \pi]\tran \in \mathbb{R}^{2N} \;.
\end{equation}
Alternative option could be to first pre-calculate a suitable trajectory by solving a boundary value problem and then use the MPC to track the trajectory. 
In our case, this option was not necessary.
We identify the predictor using the same approach as described in the previous section, differencing only in linearization of the FK model around $x_\mathrm{unst}$, instead of $x_\mathrm{stab}$.

\begin{figure}[t]
        \centering      
        \includegraphics[width=0.7\textwidth]{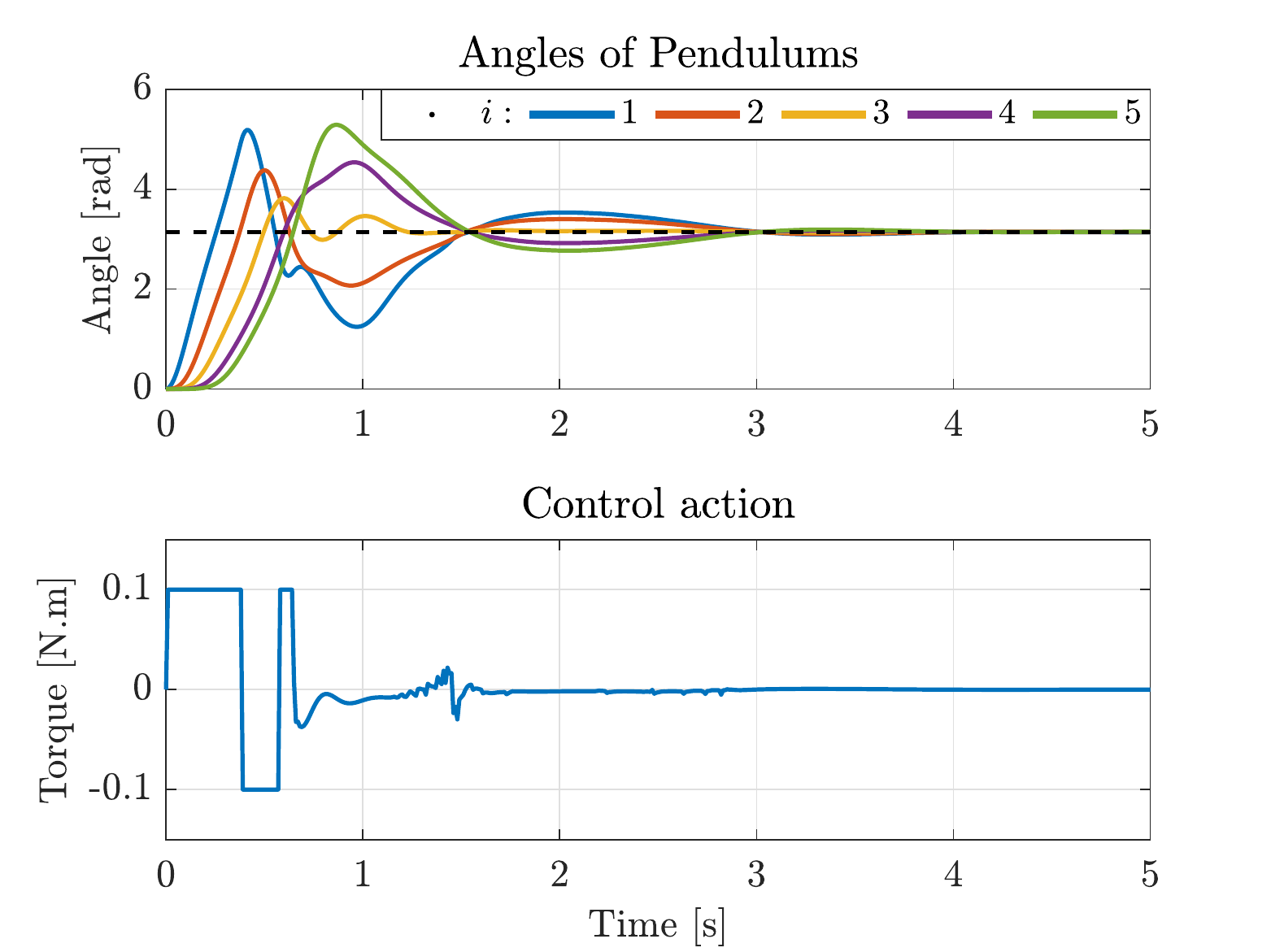}
        \caption{Controlled synchronization to the unstable equilibrium $x_\mathrm{unst}$. The dashed line represents the reference angle.}
        \label{fig:sim_ctrl_5_pend_unstable_equi}
\end{figure}

For the MPC, we selected the same cost matrices as in the previous task, thus~\eqref{eq:Q_R_matrices_equilibrium_stab} with the $Q$ matrix~\eqref{eq:Q_matrix_increasing}.
The resulting simulation is in Fig.~\ref{fig:sim_ctrl_5_pend_unstable_equi}.

\subsubsection{Periodic Trajectory}\label{sec:task_near_sync_tracking}

\begin{figure}[t]
        \centering      
        \includegraphics[width=0.7\textwidth]{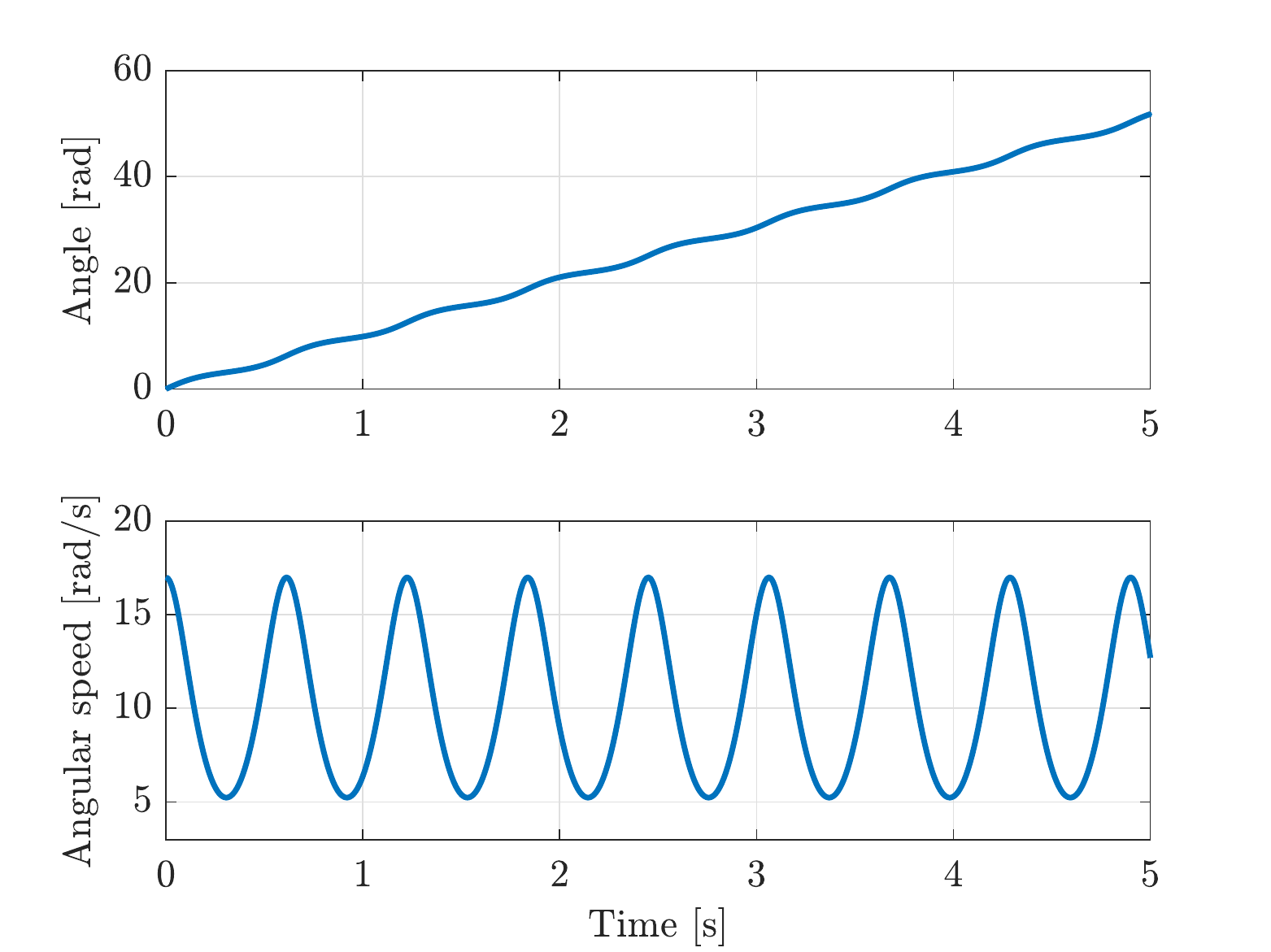}
        \caption{Periodic solution of the drift dynamics~\eqref{eq:drift_dynamics_MAS} with $\gamma = 0$ and the initial condition $x_0 = [0, 17]\tran$}
        \label{fig:ref_traj_periodic_sync}
\end{figure}

Consider now the drift dynamics of a single, \textit{virtual}, uncoupled pendulum~\eqref{eq:drift_dynamics_MAS} with $\gamma = 0$.
Additionally, consider such initial conditions, that the pendulum has enough energy to swing though the inverse position, resulting in continuous, periodic rotation.
For instance, we can generate such a trajectory with an initial condition  $x_\mr{vir}(t_0) = [0, 17]\tran$; see Fig.~\ref{fig:ref_traj_periodic_sync}.
Note, that with $\gamma > 0$, the energy of the system would naturally dissipate, so the pendulum would approach the stable equilibrium in the downward position.

Since we consider the FK model with $\gamma > 0$, and the proposed periodic solution $x_\mr{vir}(t)$ is generated with $\gamma = 0$, there is a mismatch between the dynamics of the virtual leader and the drift dynamics of the pendulums.
Therefore, the synchronization cannot be exact as the assumption stated in Sec.~\ref{sec:problem_statement} is not satisfied.
Nevertheless, on a short time horizon, the dissipation term is not dominant and, thus, we can view it only as a disturbance to the synchronization.
We then rely on the applied feedback control to mitigate the deviations from the synchronous state.

The goal is to satisfy~\eqref{eq:main_problem_statement} with $x^\star(t) = r(t)$ where $r(t)$ is the periodic solution to the drift dynamics~\eqref{eq:drift_dynamics_MAS} with $\gamma = 0$.
In other words, the goal is to synchronously rotate all the pendulums in the FK model. 
The proposed task is motivated by a potential option to decrease a nanoscale friction using a active control, see~\cite{do_synchronization_2021} for details.

To get the time-series of the reference signal $r_k = [\varphi_{\mr{vir},k}, \dot{\varphi}_{\mr{vir},k}]\tran$, we numerically solve the system~\eqref{eq:drift_dynamics_MAS} with $\gamma = 0$ and initial conditions
\begin{equation}
        [\varphi_\mr{vir}(t_0), \dot{\varphi}_{\mr{vir}(t_0)}]\tran = [0, 17]\tran \;,
\end{equation}
and integration step $\Delta t = \SI{5}{\milli\second}$.

For the identification of the predictor, we use only a proportional controller since the LQR with the linear approximation around a single point is not applicable.
At the time step $k$, we choose to control the angle of the first pendulum in the chain $\varphi_{1}$ to the reference angle $\varphi_{0}$.
Thus, the control action to gather identification data is
\begin{equation}\label{eq:PID_periodic_sync}
        u_k =  k_\mathrm{p} \left(\varphi_{\mr{vir},k} - \varphi_{1,k} \right) + v_k\;, \quad v_k \sim \mathcal{N}(0, 0.01) \;,
\end{equation}
with the gain $k_\mathrm{p} = 0.2$ and $v_k$ is again a random perturbation to get richer trajectories.  
We again gather $N_\mathrm{traj} = 100$ trajectories and identify the predictor.

\begin{figure}[t]
        \centering      
        \includegraphics[width=0.7\textwidth]{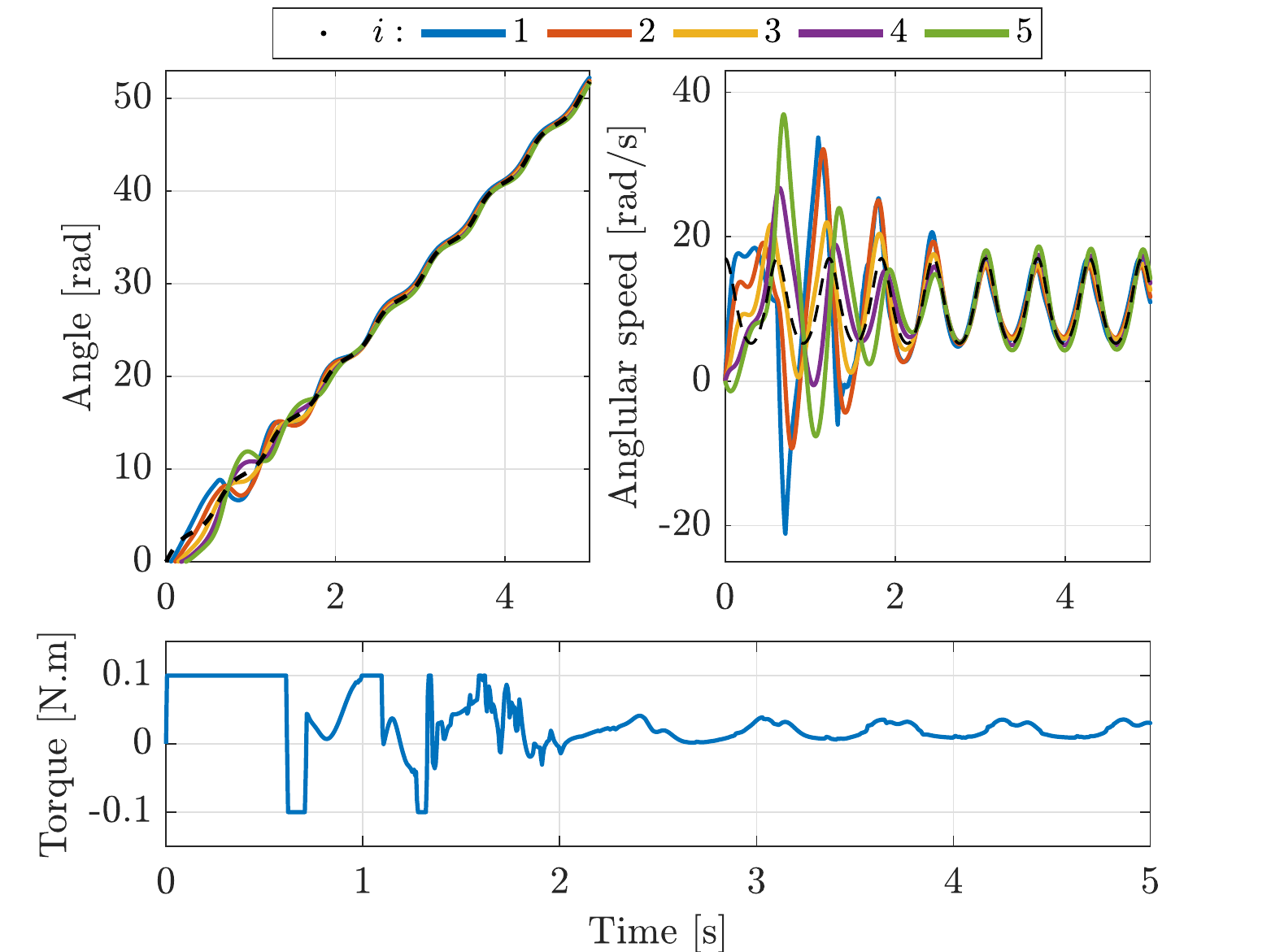}
        \caption{Controlled synchronization to a periodic trajectory $x_\mathrm{per}$ with $N = 5$ pendulums. The dashed lines represent the reference $r_k$.} 
        \label{fig:sim_ctrl_5_pend_periodic}
\end{figure}

For the MPC, we choose the same matrix $R$ as in previous tasks but for we choose different matrix $Q$.
The motivation is to penalize not only the errors $e_k$, but directly the differences between the pendulums' states.
We reflect this by penalizing the dissipated energy in a relative speeds of the pendulums
\begin{equation}\label{eq:dissipated_energy}
        D_\mr{rel} = \frac{1}{2} b \sum_{i = 1}^{N-1} \left( \dot{\varphi}_{i+1} - \dot{\varphi}_i \right)^2\;.
\end{equation}
By expanding the term $\left( \dot{\varphi}_{i+1} - \dot{\varphi}_i \right)^2$, one can see, that the dissipated energy can be included into the cost function $J$ by changing the cost matrix
\begin{equation}
        Q = I_N \otimes
        \begin{bmatrix}
                10 i^3        & 0 \\
                0               & 1
        \end{bmatrix} \; 
        +
        q_\mathrm{s}
        \frac{1}{2}b 
        \left(
        L \otimes
        \begin{bmatrix}
                0 & 0 \\
                0 & 1 \\
        \end{bmatrix}
        \right) \;,
\end{equation}
where the first term again penalizes the pendulums' states and the second term penalizes the total dissipated energy weighted with a constant ${q_\mathrm{s} > 0}$.
The resulted control is in Fig.~\ref{fig:sim_ctrl_5_pend_periodic}. 
We can see that the system approximately synchronizes after several periods.

\subsection{Hardware Experiments}\label{sec:experiments}

\begin{figure}[t]
        \centering      
        \includegraphics[width=0.5\textwidth]{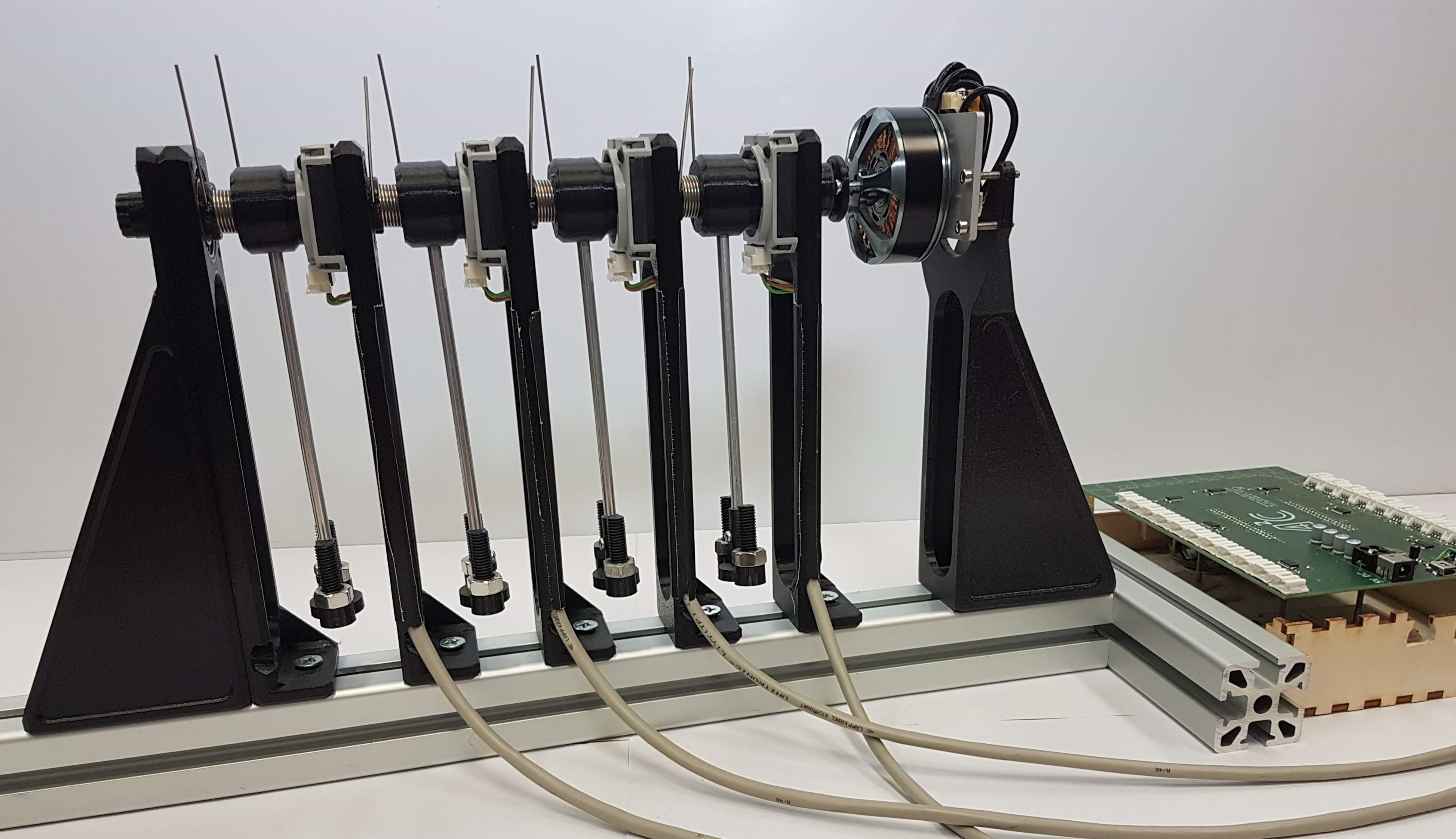}
        \caption{Mechanical FK model with $N=4$ pendulums}
        \label{fig:FK_platform}
\end{figure}

The hardware experiments were conducted on the platform depicted in Fig.~\ref{fig:FK_platform}.
The platform consists of pendulums coupled with steel torsion springs, a high-power brushless DC motor, and electronics for control and data acquisition.
The angles of pendulums are read electronically using high-resolution rotary capacitive encoders.
Further details of the platform's design are presented in~\cite{do_experimental_2022}.

\subsubsection{Results}

In contrast to simulations, we carried out the experiments only with $N = 4$ pendulums.
Experiments with four pendulums gave acceptable and repeatable results, whereas satisfying the goals with a higher number of pendulums was challenging due to the mechanical imperfection of the platform's design.

Since collecting many closed-loop trajectories for predictor's identification from the hardware platform would be time-consuming, we used the method of \textit{digital twin}.
Thus, we gathered the closed-loop trajectories, directly using the KMPC.
For all experiments, we gathered $N_\mathrm{traj} = 100$ trajectories to identify the predictor.

\begin{figure}[t]
        \centering      
        \includegraphics[width=0.7\textwidth]{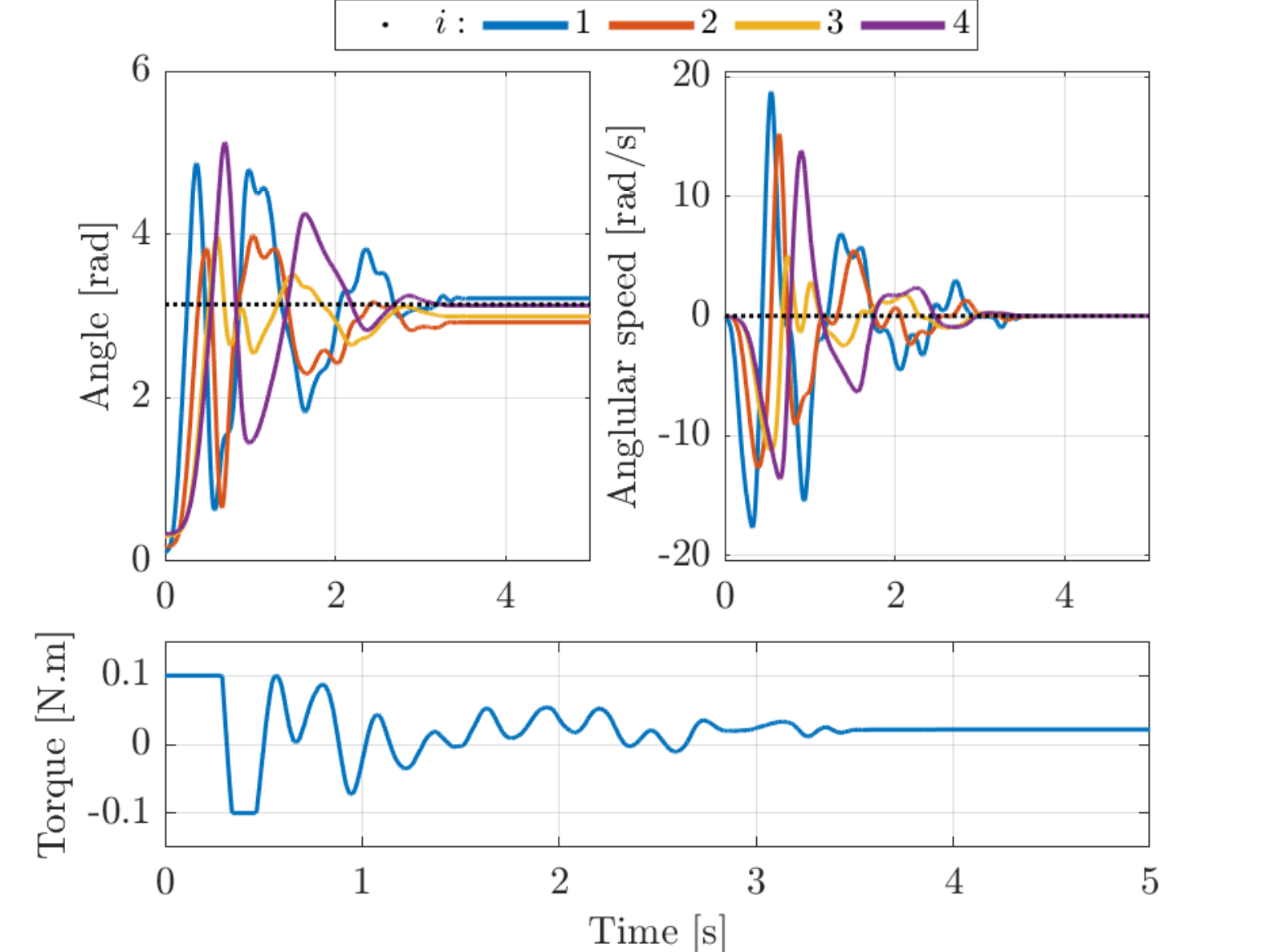}
        \caption{Hardware experiment: controlled synchronization to an unstable equilibrium $x_\mathrm{unst}$ with $N = 4$ pendulums}
        \label{fig:experiment_swing_up_4_pend}
\end{figure}

\begin{figure}[t]
        \centering      
        \includegraphics[width=0.55\textwidth]{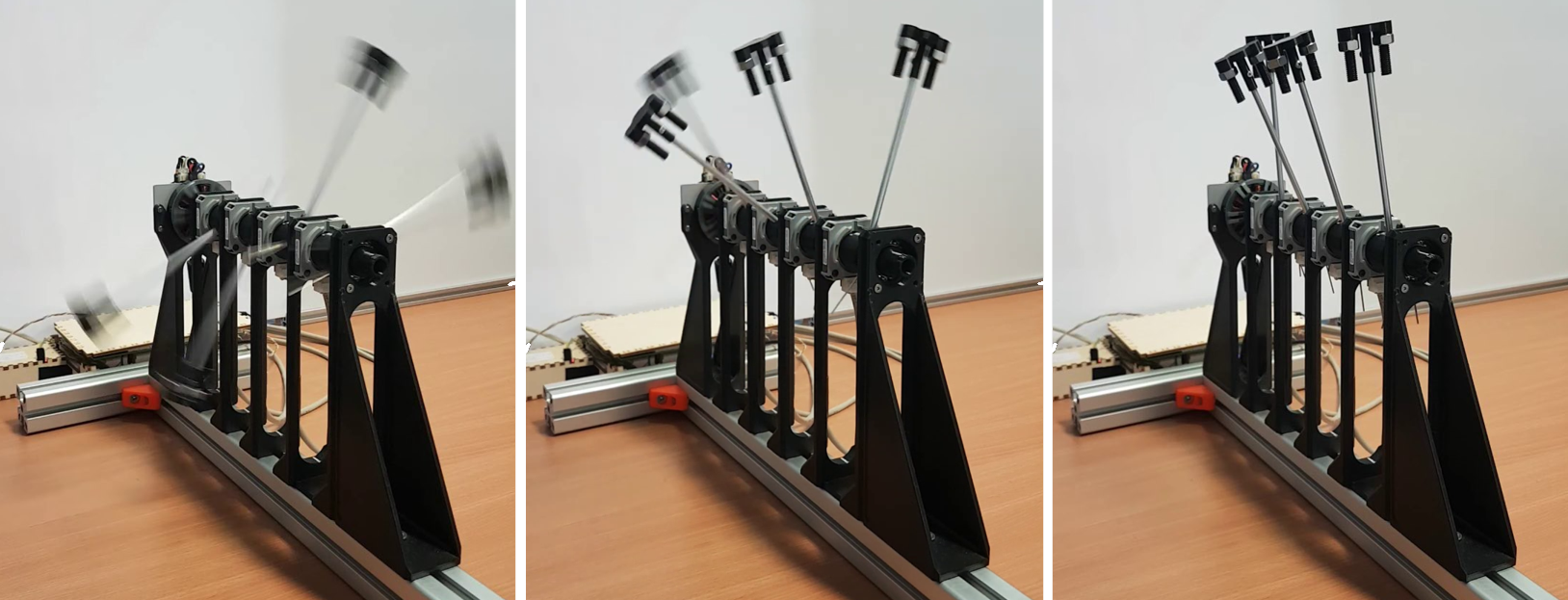}
        \caption{Snapshots of the synchronization to the unstable equilibrium with $N = 4$ pendulums}
        \label{fig:Swingup_experiment_snapshots}
\end{figure}

\begin{figure}[t]
        \centering      
        \includegraphics[width=0.7\textwidth]{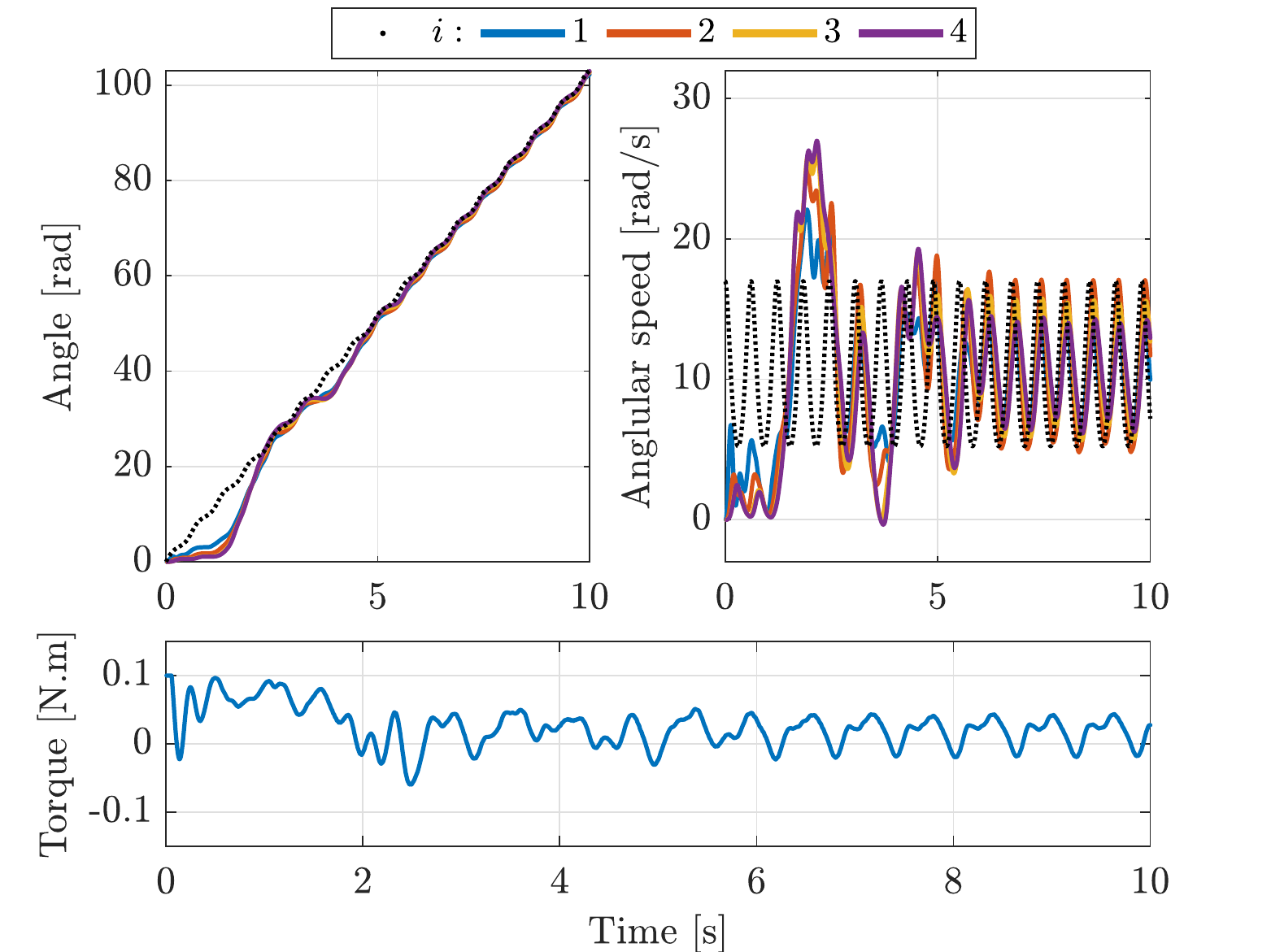}
        \caption{Hardware experiment: controlled synchronization to a periodic trajectory $x_\mathrm{per}$ with $N = 4$ pendulums. The dashed lines represent the reference $r_k$.}
        \label{fig:exper_sync_5_pends}
\end{figure}

\begin{figure}[t]
        \centering      
        \includegraphics[width=0.55\textwidth]{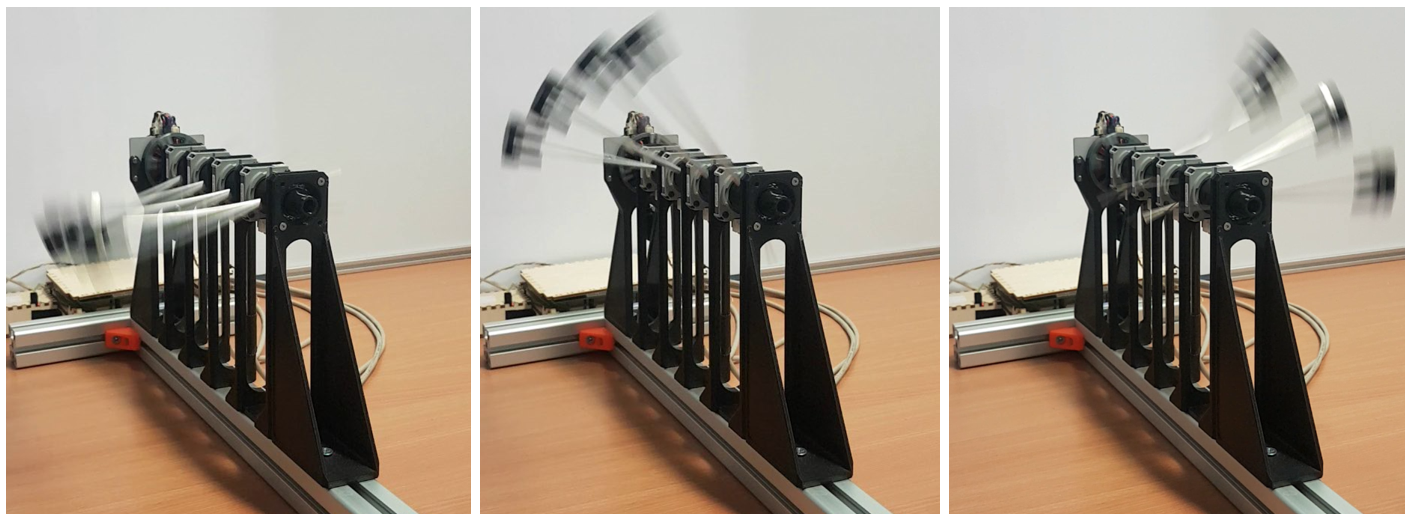}
        \caption{Snapshots of the synchronization to a periodic solution with $N = 4$ pendulums}
        \label{fig:Sync_experiment_snapshots}
\end{figure}

We showed two experiments, synchronization to the unstable equilibrium and the periodic trajectory as described in~\ref{sec:task_swing_up} and~\ref{sec:task_near_sync_tracking}.
The results are depicted in Fig.~\ref{fig:experiment_swing_up_4_pend} and~\ref{fig:exper_sync_5_pends}, with snapshots in Fig.\ref{fig:Swingup_experiment_snapshots} and~\ref{fig:Sync_experiment_snapshots}, respectively.
In comparison to the simulations, the synchronization errors were higher.
This was caused by a combination of several issues connected with hardware experiments.
For instance, the springs connecting the pendulums have a mechanical dead zone, i.e., a small deviation in a spring does not result in a proportional torque but is zero or significantly smaller.
Some problems could also be caused by the delay in the control loop, inevitably created by the hardware communication.
Nevertheless, the experimental results show that the method is robust even for implementation on the hardware platform with $\SI{5}{\milli\second}$ control period.




\section{Conclusions and Future Work}\label{sec:future_work}

In this work, we presented controlled synchronization of the FK model and successfully demonstrated the solution in simulations and on the hardware platform.
Specifically, we showed synchronization to the stable equilibrium, the unstable equilibrium (swing-up), and the periodic orbit.
We formulated the tasks as a special case of trajectory tracking and used the Koopman Model Predictive Control to solve them.
The MATLAB implementation of the simulations are available at \href{https://github.com/aa4cc/KoompanMPC-for-synchronization}{github.com/aa4cc/KoompanMPC-for-synchronization}.

The future work will aim at imposing guarantees of the controller's design.
In particular, we will focus on analyzing the identified predictor.
Another possible future direction might be to compare several sets of observables to find a better predictor for the system's behavior.
In turn, better predictions of the system's dynamics could allow us to control the FK model with more pendulums.

\section*{Appendix}
\subsection{Dense formulation}\label{sec:Appdx:dense_MPC}
The matrices and vectors of the dense formulation~\eqref{eq:MPC_formulation_dense} are
\begin{equation}
\begin{aligned}
    H &= \hat{R} + \hat{B}\tran \hat{C}\tran \hat{Q} \hat{C} \hat{B} \;,    
    \\
    \bar{{A}} &= 
        \begin{bmatrix}
            I_{Nn \times Nn} \\
            \hat{A} \\
            0_{Nn \times Nn}
        \end{bmatrix}                           \;, \quad
    \bar{{B}} =
    \begin{bmatrix}
        0_{1 \times N_\mr{p}} \\
        \hat{B} \\
        I_{N_{\mr{p}} \times N_{\mr{p}}}
    \end{bmatrix} \;,
    \\ 
    b_\mr{min} &= 
        \begin{bmatrix}
            z_\mr{min} \\
            u_\mr{min}
        \end{bmatrix}                           \;, \quad
    b_\mr{max} = 
        \begin{bmatrix}
            z_\mr{max} \\
            u_\mr{max}
        \end{bmatrix} \;,    
\end{aligned}
\end{equation}
where
\begin{equation}
    \begin{aligned}
    \hat{C} &= I_{N_\mathrm{p}} \otimes C \;,  \quad \hat{R} = I_{N_\mathrm{p}} \otimes R \;, \quad \hat{Q} = I_{N_\mr{p}}  \otimes Q \;,
    \\
    \hat{A} &= 
        \begin{bmatrix}
            A \\
            \vdots \\
            A^{N_\mathrm{p}}
        \end{bmatrix},\qquad
    \hat{B} = 
        \begin{bmatrix}
            B                   & 0         & \ldots & 0 \\
            AB                  & B         & \ldots & 0 \\
            \vdots              & \ddots    & \ddots &   \\
            A^{N_\mathrm{p}-1}B & \ldots    & AB     & B 
        \end{bmatrix} 
    \;, 
    \end{aligned}
\end{equation}


\section{Acknowledgements}
This work has been supported by the Czech Science Foundation (GACR) under contract No. 20-11626Y, by the AI Interdisciplinary Institute ANITI funding, through the
French “Investing for the Future PIA3” program under the Grant agreement n$^\circ$ ANR-19-PI3A-0004 as well as by the National Research Foundation, Prime Minister’s Office, Singapore, under its Campus for Research Excellence and Technological Enterprise (CREATE) programme.

\bibliographystyle{IEEEtran}
\bibliography{biblio_files/biblio, biblio_files/biblio_MK}

\end{document}